\documentclass[smallextended]{svjour3}
\usepackage{graphicx}
\usepackage{epsfig,amssymb}
\newcommand{\Diag}{{\mbox{Diag }}}
\newcommand{\barbx}{{\bar{\bf x}}}
\newcommand{\bxi}{{\mbox{\boldmath$\xi$}}}
\newcommand{\PPd}{P^d}
\newcommand{\sta}{{\rm sta}}
\newcommand{\vsig}{\varsigma}
\newcommand{\calP}{{\cal{P}}}
\newcommand{\calY}{{\cal{Y}}}
\newcommand{\calC}{{\cal{C}}}
\newcommand{\calZ}{{\cal{Z}}}

\newcommand{\calX}{{\cal{X}}}
\newcommand{\calS}{{\cal{S}}}
\newcommand{\calE}{{\cal{E}}}
\newcommand{\bx}{{\bf x}}
\newcommand{\bu}{{\bf u}}
\newcommand{\bz}{{\bf z}}
\newcommand{\bG}{{\bf G}}
\newcommand{\bF}{{\bf F}}
\newcommand{\half}{\frac{1}{2}}
\newcommand{\QQ}{Q}
\newcommand{\GG}{G}
\newcommand{\PP}{P}
\newcommand{\DD}{D}
\newcommand{\VV}{V}
\newcommand{\UU}{U}
\newcommand{\ZZ}{Z}
\newcommand{\BB}{B}
\newcommand{\KK}{K}
\newcommand{\WW}{W}
\newcommand{\xx}{x}
\newcommand{\HH}{H}
\newcommand{\NN}{N}
\newcommand{\bc}{{\bf c}}
\newcommand{\be}{{\bf e}}
\newcommand{\bh}{{\bf h}}
\newcommand{\by}{{\bf y}}
\newcommand{\bb}{{\bf b}}
\newcommand{\bA}{{\bf A}}
\newcommand{\bQ}{{\bf Q}}
\newcommand{\bB}{{\bf B}}
\newcommand{\bD}{{\bf D}}
\newcommand{\bH}{{\bf H}}
\newcommand{\real}{{\mathbb R}} 

\newcommand{\beps}{{\mbox{\boldmath$\epsilon$}}}
\newcommand{\brho}{\mbox{\boldmath$\rho$}}
\newcommand{\bdelta}{\mbox{\boldmath$\delta$}}
\newcommand{\bdel}{\mbox{\boldmath$\delta$}}
\newcommand{\barby}{\bar{\bf y}}
\newcommand{\bg}{{\bf g}}
\newcommand{\bff}{{\bf f}}
\newcommand{\barx}{\bar{x}}
\newcommand{\bary}{\bar{y}}
\def\la{\langle}
\def\ra{\rangle}
\newcommand{\bveps}{\mbox{\boldmath$\varepsilon$}}
\newcommand{\bvsig}{\mbox{\boldmath$\varsigma$}}
\newcommand{\btau}{\mbox{\boldmath$\tau$}}

\newcommand{\barbsig}{\bar{\mbox{\boldmath$\sigma$}}}
\newcommand{\bsig}{\mbox{\boldmath$\sigma$}}
\newcommand{\bmu}{\mbox{\boldmath$\mu$}}
\newcommand{\barbvsig}{\bar{\mbox{\boldmath$\varsigma$}}}
\newcommand{\barbtau}{\bar{\mbox{\boldmath$\tau$}}}
\newcommand{\barbmu}{\bar{\mbox{\boldmath$\mu$}}}
\newcommand{ \Lam}{{\Lambda}}
\newcommand{\calXa}{{\calX_a}}
\newcommand{ \lam}{{\lambda}}
\newcommand{\alp}{{\alpha}}

\begin{document}
\title{Global Optimal Solution to  Discrete Value Selection Problem
with Inequality Constraints}

\author{Ning Ruan        \and
        David Yang Gao
}
\institute{Ning Ruan \at
              School of Sciences, Information Technology and Engineering, \\
              University of Ballarat, Ballarat, VIC 3353, Australia. \\
              and\\
              Department of Mathematics and  Statistics, \\
              Curtin University, Perth,  WA 6845, Australia.\\
              Tel.: +613-53279942\\
              \email{n.ruan@ballarat.edu.au}           
           \and
           David Yang Gao \at
              School of Sciences, Information Technology and Engineering, \\
              University of Ballarat, Ballarat, VIC 3353, Australia. \\
              Tel.: +613-53279791\\
              \email{d.gao@ballarat.edu.au}
}

\date{Received: date / Accepted: date}

\titlerunning{Discrete Value Selection Problem}
\authorrunning{Ning Ruan and David Y. Gao}
\maketitle

\begin{abstract}
This paper presents a canonical dual method for solving a quadratic
discrete value selection problem subjected to inequality constraints.
The problem is first transformed into a problem with quadratic objective
and 0-1 integer variables. The dual problem of the 0-1 programming problem is
thus constructed
by using the canonical duality theory. Under appropriate conditions,
this dual problem
is a maximization problem of a concave function over a convex continuous space.
Numerical simulation studies, including some large
scale problems, are carried out so as to demonstrate the effectiveness and
efficiency of the method proposed.
\end{abstract}

\section{Introduction}
Many decision making problems, such as portfolio
selection, capital budgeting, production planning, resource allocation,
and computer networks,  can often be formulated as
integer programming problems. See for
examples, (Chen et al 2010; Floudas 2000; Karlof 2006).
In engineering applications, the variables of these optimization problems can
not have arbitrary values. Instead, some or all of the variables
must be selected from a list of integer or discrete values for
practical reasons. For examples, structural members may have to be
selected from selections available in standard sizes, member thicknesses may
have to be selected from the commercially available ones, the
number of bolts for a connection must be an integer, the number of
reinforcing bars in a concrete member must be an integer,etc (Huang and Arora 1997).
However, these integer programming problems
are computationally highly demanding. Nevertheless, some numerical methods are
now available.

Several review articles on
nonlinear optimization problems with discrete variables have
been recently published (Arora et al. 1994; Loh and Papalambros 1991;
Samdgren 1990; Thanedar and Vanderplaats 1994),
and some popular methods have been discussed, including
branch and bound methods,  a hybrid method that combines a
branch-and-bound method with a dynamic programming technique
(Marsten and Morin 1978), sequential linear programming,
rounding-off techniques, cutting plane techniques (Balas et al. 1993), heuristic techniques,
penalty function approach
and sequential linear programming. The relaxation method has also been proposed, leading to
second order cone programming (SOC) (Ghaddar et al. 2011). More recently, simulated
annealing (Kincaid and Padula 1990) and genetic algorithms have been discussed.

Branch and bound is perhaps the most widely known and used method for
discrete optimization problems. When applied to linear problems, this method can
be implemented in a way to yield a global minimum point;
however, for nonlinear problems there is no
such guarantee, unless the problem is convex. The branch and bound method has been used successfully
to deal with problems with discrete design variables, however, for the problem with
a large number of discrete design variables, the number of
subproblems (nodes) becomes large, making the method inefficient.

Simulated annealing (SA) is a stochastic technique to find a global minimizer. The basic idea of the method is to generate a random point and evaluate the problem functions. If the trial point is feasible and the cost function value is smaller than the current best record, the point is
accepted, and record for the best value is updated. The acceptance is based on value of the
probability density function. In computing the probability a
parameter called the temperature is used. Initially, a larger target value is selected. As the
trials progress, the target value is reduced (this is called the cooling schedule), and the process
is terminated after a fairly large number of trials. The main deficiency of the method is the unknown
rate at which the target level is to be reduced and uncertainty in the total number of trials.

Genetic algorithms (GA) belong to the category of stochastic search methods (Holland 1975). In a GA,
several design alternatives, called a population in a generation, are allowed to reproduce and
cross among themselves, with bias allocated to the most fit members of the population.
Three operators are needed to implement the
algorithm: reproduction, crossover, and mutation.
These three steps are repeated for successive generations
of the population until certain stopping criteria are satisfied.
The member in the final generation with the best fitness level is the optimum design.

The SA method and GA usually need large execution times
to find a global minimum. Although it is possible to find the best solution if temperature is
reduced slowly and enough execution time is allowed.
A drawback of SA is lack of an effective
stopping criterion. It is difficult to tell whether a global or fairly good solution has been
reached. It is also important to note that the CPU times for SA and GA can vary from one run to
the next for the same problem (Huang and Arora 1997).

Canonical duality theory provides a new and potentially useful methodology
for solving a large class of integer programming problems. It was shown
in  (Gao 2007 and  Fang et al 2008) that  the Boolean integer programming
problems are actually  equivalent to certain canonical dual problems
in continuous space  without duality gap, which can be solved
deterministically  under certain conditions.
This theory has been generalized for solving multi-integer programming
and the well-known max cut problems (see Wang et al  2008 and Wang et al 2012).
It is also shown in (Gao 2009,  Gao and Ruan 2010) that by the canonical duality
theory, the NP-hard quadratic integer programming problem can be transformed
 to a continuous unconstrained Lipschitzian
global optimization problem, which can be solved via deterministic methods
(see Gao  et al 2012).

In this paper, our goal is to solve  a general quadratic programming problem with
its decision variables taking values from discrete sets. The elements from
these discrete sets are not required to be
binary or uniformly distributed. An effective numerical method is
developed based on the
canonical duality theory (Gao 2000). The rest of the paper is organized as
follows. Section 2 presents a mathematical  statement of the problem.
Section 3 shows that this  general discrete-value quadratic programming
problem can be transformed into a 0-1 programming problem in
 higher dimensional space.
In Section 4, the canonical duality is utilized to
construct the canonical dual problem. The computational method,
which is based on
solving the canonical dual problem, is developed. Some numerical
 examples are illustrated to demonstrate the effectiveness and
efficiency of the proposed method.The paper is ended with some  concluding remarks.

\section{Discrete Programming Problem}
The discrete programming problem  to be addressed is given below:
\begin{eqnarray}\label{pri}
&(\calP_{a})\;\;&{\rm Minimize}\;\;P(\bx)=\half \bx^T \QQ \bx -\bc^T \bx\\
&& {\rm subject\;to}\;\;\bg(\bx)=\bA \bx - \bb \le 0,\\
&& \;\;\;\;\;\;\;\;\;\;\;\;\;\;\;\;\;\bx=[x_1, x_2, \cdots, x_n]^T, \;x_i \in U_i,\;i=1,\cdots,n, \nonumber
\end{eqnarray}
where
$\QQ = \{ q_{ij} \}  \in \real^{n\times n}$ is an $n\times n$ positive semi-definite
symmetric matrix,
 $\bA =\{ a_{ij} \} \in \real^{m \times n}$
is an $m\times n$ matrix with $rank(\bA)=m<n$, $\bc =[c_1,\cdots,c_n]^T  \in \real^n$
and $\bb= [b_1, \cdots, b_m]^T \in \real^m$ are given vectors.
Here, for each $i=1, \cdots, n$,
\begin{eqnarray*}
U_i=\{ u_{i,1}, \cdots, u_{i,k_i}\},
\end{eqnarray*}
where, $u_{i,j}, j=1, \cdots, K_i$, are given
real numbers. Let $\KK=\sum_{i=1}^n K_i$.
\section{Equivalent Transformation}
Let us introduce the following transformation,
\begin{eqnarray}\label{01transfer}
x_i=\sum_{j=1}^{K_i} u_{i,j} y_{i,j},\; i=1, \cdots, n,
\end{eqnarray}
where, for each $i=1, \cdots,n $, $u_{i,j} \in U_i,\; j=1, \cdots, K_i$. Then,  the
discrete programming problem $(\calP_{a})$ can be written as the following
0-1 programming problem:
\begin{eqnarray}
&(\calP_{b})\;\;&{\rm Minimize}\;\;P(\by)=\half \by^T B \by -\bh^T \by\\
&& {\rm subject\; to}\;\;\bg(\by) = \DD \by - \bb \le 0,\label{p2con1}\\
&& \;\;\;\;\;\;\;\;\;\;\;\;\;\;\;\;\;\sum_{j=1}^{K_i} y_{ij} - 1 = 0,\; i=1, \cdots, n,\label{p2con2}\\
&& \;\;\;\;\;\;\;\;\;\;\;\;\;\;\;\;\;y_{i,j} \in \{0, 1 \}, ~i=1,\ldots,n;\;
j=1,\cdots, K_i,\label{p2con3}
\end{eqnarray}
where
\begin{eqnarray*}
\by=[y_{1,1}, \cdots, y_{1, K_1},\cdots, y_{n,1}, \cdots,y_{n, K_n}]^T
\in \real^{K},
\end{eqnarray*}
\begin{eqnarray*}
\bh=[c_1 u_{1,1},\cdots, c_1 u_{1, K_1}, \cdots,
c_n u_{n,1}, \cdots, c_n u_{n, K_n}]^T
\in \real^{K},
\end{eqnarray*}
\begin{eqnarray*}
\BB=
\left[
\begin{array}{ccccc}
q_{1,1} u_{1,1}^2 & \cdots& q_{1,1} u_{1,1} u_{1, K_1}&
\cdots& q_{1,n} u_{1,1}u_{n, K_n} \\
\vdots&\ddots&\vdots&\ddots&\vdots\\
q_{1,1}u_{1, K_1} u_{1,1}&\cdots&
q_{1,1}u_{1, K_1}^2 & \cdots& \cdots\\
\vdots&\ddots&\vdots&\ddots&\vdots\\
q_{n,1} u_{n, K_n} u_{1,1}& \cdots & \cdots & \cdots & q_{n,n}u_{n, K_n}^2
\end{array}
\right]
\;\in \real^{\KK \times \KK},
\end{eqnarray*}
\begin{eqnarray*}
\DD=
\left[
\begin{array}{ccccc}
a_{1,1} u_{1,1} & \cdots & a_{1,1} u_{1, K_1}&\cdots& a_{1,n} u_{n, K_n}\\
\vdots&\ddots&\vdots&\ddots&\vdots\\
a_{m,1} u_{1,1} & \cdots & a_{m,1} u_{1, K_1}&\cdots& a_{m,n} u_{n, K_n}
\end{array}
\right]
\;\in \real^{m \times \KK}.
\end{eqnarray*}
\begin{theorem}
Problem $(\calP_{b})$ is equivalent to Problem $(\calP_{a})$.
\end{theorem}
\noindent {\em Proof}.
For any $i=1,2, \cdots, n$, it is clear that constraints
(\ref{p2con2}) and (\ref{p2con3}) are equivalent to the existence of only
one $j \in \{ 1, \cdots, K_i\}$, such that $y_{i,j}=1$
while $y_{i,j}=0$ for all other $j$. Thus, from the definition of $ \by$,
the conclusion
follows readily.
\hfill $\Box$\\
Let
\begin{eqnarray*}
\HH =
\left[
\begin{array}{cccccccccc}
1  & \cdots & 1 & 0 &\cdots& 0 &\cdots& 0 &\cdots &0\\
0  & \cdots & 0 & 1&\cdots& 1 &\cdots& 0 &\cdots &0\\
\vdots&\ddots&\vdots&\vdots&\ddots&\vdots&\ddots&\vdots&\ddots&\vdots\\
0  & \cdots & 0 & 0 &\cdots& 0 &\cdots& 1 &\cdots &1
\end{array}
\right]\in \real^{n \times \KK}
\end{eqnarray*}
and, for any integer $N$, let
\begin{eqnarray*}
\be_{\NN}=[1, \cdots, 1,\cdots, 1, \cdots,1]^T
\in \real^{\NN}.
\end{eqnarray*}
We consider the following quadratic programming problem:
\begin{eqnarray}
&(\calP)\;\;&{\rm Minimize}\;\;P(\by)=\half \by^T B \by -\bh^T \by\\
&& {\rm subject\; to}\;\;\bg(\by) = \DD\by -\bb \le 0,\\
&& \;\;\;\;\;\;\;\;\;\;\;\;\;\;\;\;\;\HH \by - \be_n = 0,\\
&& \;\;\;\;\;\;\;\;\;\;\;\;\;\;\;\;\;\by\circ(\by-\be_K)\le 0,
\end{eqnarray}
where the notation
${\bf s}\circ{\bf t}:=[s_1 t_1,s_2 t_2,\ldots,s_K t_K]^T$
denotes  the Hadamard product for any two vectors
${\bf s}, {\bf t}  \in \real^K$.

\section{Canonical duality theory: A brief review}
The basic idea of the canonical duality theory can be   demonstrated
by solving the following
general nonconvex problem (the primal problem $(\calP)$ in short)
\begin{eqnarray}
(\calP): \; \min_{ \xx \in \calX_a}
\left\{ \PP(\xx) = \half   \la \bx ,  \bA \xx  \ra  -
\la  \xx ,  \bff \ra  +  \WW( \xx) \right\},
\end{eqnarray}
where
$ \bA  \in \real^{n\times n} $ is  a given symmetric indefinite matrix,
$\bff \in \real^n$ is a given vector,  $\la \xx, \xx^* \ra $
denotes the bilinear form
between $\xx$ and its dual variable $\xx^*$,
$\calX_a \subset \real^n$ is a
given  feasible space,
and $W: \calX_a \rightarrow \real \cup \infty$ is
 a general nonconvex objective function.
The mathematical definition of the
objectivity for general functions is
given in (Gao, 2000).

The  \textbf{key step} in  the canonical dual  transformation
is to choose a nonlinear operator,
\begin{eqnarray}
\bveps = \Lam (\bx):\calX_a \rightarrow \calE_a \subset \real^p
\end{eqnarray}
and a {\em canonical function} $\VV: \calE_a \rightarrow \real$
such that the nonconvex objective function $\WW( \bx)$
can be recast by adopting a canonical form
$\WW( \bx) = \VV(\Lam(\bx))$.
Thus, the primal problem $(\calP)$  can be
written in the following canonical form:
\begin{eqnarray}
(\calP): \;  \min_{\bx \in \calX_a}
\left\{ \PP(\bx) =   \VV(\Lam(\bx)) - \UU(\bx)\right\},
\label{eq-canform}
\end{eqnarray}
where $\UU(\bx) =  \la \bx, \bff \ra - \half \la \bx , \bA \bx \ra$.
By the definition introduced in
(Gao 2000), a differentiable function $\VV(\bveps)$ is said to be
a \textit{canonical function}  on its domain  $\calE_a$ if the
duality mapping $\bvsig = \nabla \VV(\bveps)$ from $\calE_a$ to its range
$ \calS_a \subset \real^p $
is invertible. Let   $\la \bveps ;  \bvsig \ra $
denote  the bilinear form on $\real^p$.
Thus, for the given canonical function $\VV(\bveps)$,
its Legendre conjugate
$\VV^*(\bvsig)$ can be defined uniquely by the Legendre transformation
\begin{eqnarray}
\VV^*(\bvsig)  = \sta \{ \la \bveps ;  \bvsig \ra - \VV(\bveps ) \; | \; \; \bveps \in \calE_a \},
\end{eqnarray}
where  the notation $\sta \{ g(\bveps) | \; \bveps \in \calE_a\}$
stands for finding stationary point of $g(\bveps)$ on $\calE_a$.
It is easy to prove that
the following canonical  duality relations hold on $ \calE_a \times \calS_a$:

\begin{eqnarray}
\bvsig =\nabla \VV(\bveps) \; \Leftrightarrow \;  \bveps = \nabla
\VV^*(\bvsig) \; \Leftrightarrow
\VV(\bveps) + \VV^*(\bvsig) =
\la \bveps ; \bvsig \ra  . \label{eq-candual}
\end{eqnarray}
By this one-to-one canonical duality,   the nonconvex term
$W( \bx)=\VV(\Lam(\bx))$ in the problem $(\calP)$ can be replaced by
$\la \Lam(\bx) ; \bvsig \ra -
\VV^*(\bvsig)$ such that the nonconvex
function $\PP(\bx)$ is reformulated  as
the so-called Gao and Strang total complementary function (Gao 2000):
\begin{eqnarray}
\Xi(\bx, \bvsig) = \la  \Lam(\bx) ; \bvsig \ra    - \VV^*(\bvsig) -  \UU(\bx).
\label{eq:xi}
\end{eqnarray}
By using this total complementary function,
the canonical dual function $\PP^d (\bvsig)$
can be obtained   as
\begin{eqnarray}
\PP^d(\bvsig) &=& \sta \{ \Xi(\bx, \bvsig) \; | \; \bx \in \calX_a   \} \nonumber \\
&=& \UU^\Lam(\bvsig) - \VV^*(\bvsig),
\end{eqnarray}
where $\UU^\Lam(\bx)$ is defined by
\begin{eqnarray}
\UU^\Lam(\bvsig) =  \sta \{ \la  \Lam(\bx) ;  \bvsig \ra - \UU(\bx) \;
| \;\; \bx \in \calX_a \}.
\end{eqnarray}
In many applications, the geometrically nonlinear operator
$\Lam(\bx)$ is usually  quadratic function
\begin{eqnarray}
\Lam(\bx)=\half \la \bx, D_k \bx \ra +\la \bx, \bb_k \ra,
\end{eqnarray}
where $D_k \in \real^{n \times n}$ and
$\bb_k\in \real^n  (k = 1, \cdots, p)$. Let $\bvsig= [\vsig_1,\cdots, \vsig_p]^T$.
In this case, the canonical dual function can be written in the following form:
\begin{eqnarray}
\PPd(\bvsig)=-\half \la \bF(\bvsig), \bG(\bvsig) \bF(\bvsig) \ra -V^{\ast}(\bvsig),
\end{eqnarray}
where $\bG(\bvsig) = \bA+\sum_{k=1}^p \bvsig_k D_k$, and
$\bF(\bvsig)=\bff-\sum_{k=1}^p \vsig_k \bb_k$.

Let $\calS^+_a = \{\bvsig\in \real^p|\; G(\bvsig) \succeq 0\}$.
Therefore, the canonical dual problem can be proposed as
\begin{eqnarray}
(\calP^d): \;\;   \max \{ \PPd(\bvsig) | \;\; \bvsig \in \calS^+_a\} \vspace{-.3cm}.
\end{eqnarray}
which is a concave maximization problem over a convex set
$\calS^+_a \subset \real^p$.

\begin{theorem}[Gao 2000]\label{duality}
Problem $(\calP^d)$ is canonically dual to $(\calP)$ in the sense that
if $\barbvsig$ is a critical point of $\Pi^d(\bvsig)$, then
\begin{eqnarray}\label{eq-anasol}
\barbx = \bG^{\dagger}(\barbvsig) \btau(\barbvsig)
\end{eqnarray}
is a critical point of $\Pi(\bx)$ and
\begin{eqnarray}
\PP(\barbx) =  \Xi(\barbx, \barbvsig) = \PP^d(\barbvsig). \label{eq-p=d}
\end{eqnarray}
If $\barbvsig$ is a solution to $(\calP^d)$, then $\barbx  $
is a global minimizer of $(\calP) $ and
\begin{eqnarray}
\min_{\bx \in \calXa } \PP(\bx) =  \Xi(\barbx, \barbvsig)=
\max_{\bvsig \in \calS_c^+} \PP^d(\bvsig)\label{trisaddle}.
\end{eqnarray}
Conversely, if $\barbx$ is a solution to $(\calP)$, it must be in the form of
(\ref{eq-anasol}) for critical solution $\barbvsig$ of $\Pi^d(\bvsig)$.
\end{theorem}

\begin{figure}[!t]
\centering
\includegraphics[width=2.5in]{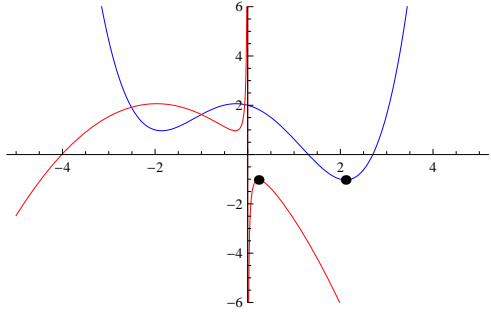}
\caption{Graphs of the primal function $\Pi(\bx)$  (blue)
and its canonical dual function $\Pi^d(\vsig)$ (red).}
\label{onedim}
\end{figure}

To help explain the theory, we   consider a simple nonconvex optimization in
$\real^n$:
\begin{eqnarray}
\min \Pi(\bx)=\half \alpha (\half\|\bx\|^2-\lam)^2-\bx^T \bff, \; \forall \bx \in \real^n,
\end{eqnarray}
where $\alp, \lam > 0$ are given parameters.
The criticality condition $\nabla P(\bx)=0$ leads to a nonlinear algebraic
equation system in $\real^n$
\begin{eqnarray}
\alpha (\half \|\bx\|^2-\lam)\bx =\bff.
\end{eqnarray}
Clearly, to solve this n-dimensional nonlinear algebraic equation directly is  difficult.
Also traditional convex optimization theory
 can't be used to identify global minimizer.  However, by the
canonical dual transformation, this problem can be solved.
To do so, we let  $\bxi=\Lam(u)=\half\|\bx\|^2-\lam \in \real$. Then,
the nonconvex function $W(\bx) = \half \alpha(\half \| \bx \|^2 -\lam)^2$
can be written in canonical form
$V(\bxi) = \half \alpha \bxi^2$.
Its Legendre conjugate is  given by
$V^{\ast}(\vsig)=\half \alpha^{-1}\vsig^2$, which is strictly convex.
Thus,
the total  complementary function for this nonconvex optimization
problem is
\begin{eqnarray}
\Xi(\bx,\vsig)=(\half \|\bx\|^2 - \lam) \vsig-\half
\alpha^{-1}\vsig^2 - \bx^T \bff.
\end{eqnarray}
For a fixed $\vsig \in \real$, the  criticality condition
$\nabla_{\bx} \Xi(\bx)=0$ leads to
\begin{eqnarray}\label{balance}
\vsig \bx-\bff=0.
\end{eqnarray}
For each  $\vsig \neq 0  $,
the  equation
(\ref{balance}) gives $\bx=\bff/\vsig$ in vector form. Substituting this into the
total complementary function $\Xi$,
the canonical dual function can be easily obtained as
\begin{eqnarray}
\Pi^d(\vsig)&=&\{\Xi(\bx,\vsig)| \nabla_{\bx} \Xi(\bx,\vsig)
=0\}\nonumber\\
&=& -\frac{\bff^T \bff}{2 \vsig}-\half \alpha^{-1} \vsig^2
-\lam \vsig, \;\;\; \forall \vsig\neq 0.
\end{eqnarray}
The critical point of this canonical function is obtained
by solving the following dual algebraic  equation
\begin{eqnarray}
(\alpha^{-1} \vsig+\lam)\vsig^2=\half \bff^T \bff.
\end{eqnarray}
For any given parameters $\alpha$, $\lam$ and the vector $\bff\in \real^n$,
this cubic algebraic equation has at most three roots
satisfying $\vsig_1 \ge 0 \ge \vsig_2\ge \vsig_3$,
and each of these roots leads to a critical point of the nonconvex function
$P(\bx)$, i.e., $\bx_i=\bff/\vsig_i$, $i=1,2,3$. By
the fact that $\vsig_ 1 \in \calS^+_a = \{ \vsig \in \real\; |\; \vsig > 0 \}$,
then Theorem 1 tells us that $\bx_1$ is
a global minimizer of $\Pi(\bx)$.
Consider one dimension problem with $\alpha= 1$, $\lam=2$, $f= \half$,
the primal function and canonical dual function
are shown in Fig. \ref{onedim}, where,  $x_1= 2.11491$ is global minimizer
of $P(\bx)$,
$\vsig_1=0.236417$ is global maximizer of $\Pi^d(\bvsig)$, and
$\Pi(x_1)=-1.02951=\Pi^d(\vsig_1)$ (See the two black dots).

The canonical duality theory was original developed from general
nonconvex systems. The canonical dual transformation can be used to convert
a nonconvex problem into a canonical dual problem without duality gap, while the
classical dual approaches may suffer from having a potential gap
(Rockafellar 1987). The complementary-dual principle provides
a unified form of analytical solutions to general
nonconvex problems  in either continuous or discrete  systems.
The canonical duality theory
has shown its potential for various classes of challenging problems.
A comprehensive review of the canonical duality theory and its applications
can be found in (Gao and Ruan 2008; Gao and Ruan 2010; Gao et al. 2012; Gao et al. 2009;
Gao et al. 2010; Ruan et al. 2010).
\section{Canonical Dual Problem}
Now we apply the canonical duality theory to integer
programming problem  presented in Section 2.

Let
\begin{eqnarray*}
U(\by) = - P(\by) = \bh^T \by-\half \by^T \BB \by,
\end{eqnarray*}
and define
\begin{eqnarray*}
\bz &=&\Lambda(\by)= [(\DD \by-\bb)^T, (\HH \by-\be_{n})^T, (\by\circ (\by-\be_K))^T]^T\\
&=& [(\beps)^T,(\bdelta)^T,(\brho)^T]^T
\in \real^{m + n + K},
\end{eqnarray*}
where $\Lambda $ is the so-called geometric operator.
Let
\begin{eqnarray*}
\WW(\bz) =
\left\{
\begin{array}{ll}
0 & \mbox{ if } \beps \le 0, \bdelta = 0, \brho \le 0,  \\
+ \infty & \mbox{ otherwise}.
\end{array}
\right.
\end{eqnarray*}
Let $\bz^*= [(\bsig)^T,(\btau)^T,(\bmu)^T]^T
\in \real^{m + n + \KK}$
be the canonical dual variables corresponding to  those from
the set
$\ZZ=\{(\beps, \bdelta, \brho): \beps \le 0, \bdelta = 0, \brho \le 0 \}.$
Then, the Fenchel super-conjugate of the function $\WW(\bz)$  is defined by
\begin{eqnarray}
\WW^{\sharp}(\bz^*) &=& \sup\{\bz^T \bz^{*}-W(\bz):\; \bz \in \ZZ \}\nonumber\\
&=& \left\{ \begin{array}{ll}
0 & \mbox{ if } \bsig \ge 0, \bmu \ge 0,  \\
+ \infty & \mbox{ otherwise}.
\end{array} \right.
\end{eqnarray}
Let
\begin{eqnarray}\label{defg}
G(\bmu)=\BB +2 \Diag (\bmu),
\end{eqnarray}
and
\begin{eqnarray}\label{deff}
F(\bsig, \btau, \bmu) = \bh - \DD^T \bsig -\HH^T \btau + \bmu.
\end{eqnarray}
Then, the total complementary function can be  obtained as:
\begin{eqnarray*}
\Xi(\by,\bsig,\btau,\bmu) &=&
\langle \Lambda (\by), \bz^*\rangle - W^{\sharp} (\bz^*)-U(\by)\nonumber\\
&=& \half \by^T \BB \by -\bh^T \by + \bsig^T (\DD \by -\bb)\nonumber\\
&& +\btau^T (\HH \by - \be_n) +\bmu^T (\by \circ (\by -\be_{\KK}))\nonumber\\
&=& \half \by^T \BB \by +\half \by^T (2 \Diag(\bmu)) \by - \bh^T \by\nonumber\\
&& (\DD^T \bsig)^T \by - \bmu^T \by + (\HH^T \btau)^T \by
-\bsig^T \bb -\btau^T \be_n\nonumber\\
&=& \half \by^T G(\bmu) \by - F^T(\bsig,\btau,\bmu) \by
-\bsig^T \bb-\btau^T \be_{n}.
\end{eqnarray*}
The critical condition $\nabla_{\by} \Xi(\by,\bsig,\btau,\bmu) = 0$ leads to
\begin{eqnarray}
\by= G^\dagger(\bmu)F(\bsig,\btau,\bmu),
\end{eqnarray}
where $ \GG ^\dagger(\bmu) $ denotes the Moore-Penrose generalized inverse
of $\GG(\bmu)=(\BB +2 \Diag (\bmu))$.

The canonical dual problem can be stated as follows:
\begin{eqnarray*}
&(\calP^{d})\;\;&{\rm Maximize}\;\;P^d(\bsig, \btau,\bmu)
=-\half F(\bsig, \btau, \bmu)^T G^\dagger(\bmu) F(\bsig, \btau, \bmu)
-\bsig^T \bb-\btau^T \be_{n} \nonumber\\
&& {\rm subject \; to}\;\; \bsig \ge 0, \bmu >0,\\
&&\;\;\;\;\;\;\;\;\;\;\;\;\;\;\;\;\;\bsig \in \real^m, \btau \in \real^n,
\bmu \in \real^{\KK}.\nonumber
\end{eqnarray*}

\begin{theorem}[Complementary-Dual Principle]\label{thm-pdg}
Problem $(\calP^{d})$ is a canonically dual to Problem  $(\calP)$
in the sense that if  $(\barbsig, \barbtau, \barbmu)  $ is  a
KKT solution of Problem $(\calP^{d})$, then the vector
\begin{eqnarray}\label{pridualrelation}
\barby (\barbsig, \barbtau, \barbmu)=
\GG ^\dagger(\barbmu)  F(\barbsig, \barbtau, \barbmu)
\end{eqnarray}
is a KKT solution of Problem $(\calP)$ and
\begin{eqnarray*}
\PP(\barby) = \PPd(\barbsig, \barbtau, \barbmu). \label{p=pdg}
\end{eqnarray*}

Moreover, if $(\barbsig, \barbtau, \barbmu)$
is a critical point of Problem $(\calP^{d})$ and $\barbmu > 0 $, then
$\barby $ is a critical point of  Problem $(\calP)$.
\end{theorem}

\noindent {\em Proof}.
By introducing the Lagrange multiplier vector $ \beps \le 0 \in \real^m$,
$\bdel \in \real^n$, and
$\brho \le 0 \in \real^{\KK}$,
the Lagrangian function associated with the dual function
$P^d(\bsig, \btau,\bmu) $ becomes
\begin{eqnarray*}
L(\bsig, \btau,\bmu, \beps, \bdel, \brho)
= P^d(\bsig, \btau,\bmu)- \beps^T \bsig + \bdel^T \btau - \brho^T \bmu.
\end{eqnarray*}
Then, the KKT conditions of the dual problem become
\begin{eqnarray*}
\frac{\partial L(\bsig, \btau,\bmu, \beps, \bdel, \brho)}{\partial \bsig}
&=& \DD \by -\bb - \beps = 0,\label{parsig}\\
\frac{\partial L(\bsig, \btau,\bmu, \beps, \bdel, \brho)}{\partial \btau}
&=& \HH \by-\be_n +\bdel = 0,\label{partau}\\
\frac{\partial L(\bsig, \btau,\bmu, \beps, \bdel, \brho)}{\partial \bmu}
&=& \by \circ (\by-\be_K)-\brho =0,\label{parmu}\\
\bsig \ge 0, \beps \le 0, \bsig^T \beps &=& 0,\label{compsig}\\
\bmu \ge 0, \brho \le 0, \bmu^T \brho &=& 0.\label{comptau}
\end{eqnarray*}
They can be written as:
\begin{eqnarray}
\DD \by \le \bb,\label{prifeasi1}\\
\HH \by - \be_n= 0,\label{prifeasi2}\\
\by(\by-\be_K) \le 0,\label{prifeasi3}\\
\bsig \ge 0, \bsig^T( \DD\by - \bb)=0,\label{dualfeasi1}\\
\bmu \ge 0, \bmu^T(\by \circ (\by -\be_{\KK}))=0.\label{dualfeasi2}
\end{eqnarray}
Specifically, if $\bmu >0$, the complementary condition leads to
$\by \circ(\by-\be_{\KK})=0$.
This proves that
if  $(\barbsig, \barbtau, \barbmu)  $ is  a
KKT solution of $(\calP^{d})$, then
(\ref{prifeasi1})-(\ref{prifeasi3}) is the so-called  primal feasibility
condition, while  (\ref{dualfeasi1})-(\ref{dualfeasi2}) is
the so-called dual feasibility condition and complementary slackness condition.
Therefore, the vector
\begin{eqnarray*}
\barby (\barbsig, \barbtau, \barbmu)=
\GG ^\dagger(\barbmu)  F(\barbsig, \barbtau, \barbmu)
\end{eqnarray*}
is a KKT solution of Problem $(\calP)$.

Again, by the complementary condition and
(\ref{pridualrelation}), we have
\begin{eqnarray*}
P^d(\bsig, \btau,\bmu)
& = & -\half F(\bsig, \btau, \bmu)^T G(\btau)^\dagger F(\bsig, \btau, \bmu)
-\bsig^T \bb -\btau^T \be_n\nonumber\\
&=& \half \by^T \BB \by - \bh^T \by +
\bsig^T (\DD \by -\bb)
+\btau (\HH \by -\be_n)
+\bmu(\by \circ (\by-\be_K))\nonumber\\
&=& \half \by^T \BB \by - \bh^T \by
= P(\by).
\end{eqnarray*}
\hfill $\Box$\\

To continue, let the feasible space $\calY$ of problem $(\calP)$ and
the dual feasible space $\calZ$ be defined by
\begin{eqnarray*}
\calY=\{\by \in \real^{K}
:\DD \by \le b, \HH \by =\be_n, \by \circ(\by -\be_K) \le 0\}
\end{eqnarray*}
and
\begin{eqnarray*}
\calZ = \{(\bsig, \btau, \bmu) \in \real^m \times
\real^n \times  \real^K:  \bsig \ge 0, \bmu >0,
F(\barbsig, \barbtau, \barbmu)  \in \calC_{ol} (\GG(\bmu))\},
\end{eqnarray*}
respectively, where $\calC_{ol} (\GG(\bmu))$ denotes the linear space spanned by the
columns of $\GG(\bmu)$.

We introduce a subset of the dual feasible space:
\begin{eqnarray}
\calZ_a^+: = \{(\bsig, \btau, \bmu) \in \calZ: G(\bmu) \succ 0 \}.
\end{eqnarray}
We have the following theorem.

\begin{theorem}\label{pdthm}
Assume that $(\barbsig,\barbtau, \barbmu)$ is a
critical point of $P^d(\bsig,\btau,\bmu)$ and
$\barby=G^\dagger(\barbmu) F(\barbsig,\barbtau,\barbmu)$.
If $(\barbsig,\barbtau, \barbmu)\in \calZ_a^+$, then
$\barby$ is a global minimizer of $P(\by)$  and
$(\barbsig,\barbtau, \barbmu)$ is a global maximizer of
$P^d(\bsig,\btau,\bmu)$ with
\begin{equation}\label{eq21}
P(\barby)=\min_{\by \in \calY}P(\by)=\max_{(\bsig,\btau,\bmu)\in
\calZ_a^+}P^d(\bsig,\btau,\bmu)
= P^d(\barbsig,\barbtau,\barbmu)
\end{equation}
\end{theorem}
\noindent \emph{Proof}

The canonical dual function
$P^d(\bsig,\btau,\bmu)$ is concave on $\calZ_a^+$.
Therefore, a critical point $(\barbsig,\barbtau,\barbmu)\in
\calZ_a^+$ must be a global maximizer of
$P^d(\bsig,\btau,\bmu)$ on $\calZ_a^+$. For any given
$(\bsig,\btau,\bmu) \in \calZ_a^+$, the complementary
function $\Xi(\by,\bsig,\btau,\bmu)$ is convex in $\by $
and concave in $(\bsig,\btau,\bmu)$, the critical
point $(\barby,\barbsig,\barbtau,\barbmu)$
is a saddle point of the complementary function. More specifically, we have
\begin{eqnarray*}
P^d(\barbsig,\barbtau,\barbmu)
& = & \max_{(\bsig,\btau,\bmu)\in \calZ_a^+} P^d(\bsig,\btau,\bmu)\nonumber\\
& = & \max_{(\bsig,\btau,\bmu)\in \calZ_a^+}
\min_{\by \in \calY}\Xi(\by,\bsig,\btau,\bmu)
=\min_{\by \in \calY}
\max_{(\bsig,\btau,\bmu)\in \calZ_a^+}\Xi(\by,\bsig,\btau,\bmu)\nonumber\\
&=& \min_{\by \in \calY} \max_{(\bsig,\btau,\bmu) \in \calZ_a^+}
\{\half \by^T G(\bmu) \by -(\bh-\DD^T \bsig-\HH^T \btau +\bmu)^T \by \\
&&-\bsig^T \bb-\btau^T \be_n \}\nonumber\\
&=& \min_{\by \in \calY} \max_{(\bsig,\btau,\bmu)\in \calZ_a^+}
\{ \half \by^T \BB \by -\bh^T \by
+\bsig^T(\DD \by - \bb)\\
&&+\btau^T (\HH \by -\be_n)
+\bmu^T \by \circ (\by - \be_{\KK})\}\nonumber\\
&=& \min_{\by \in \calY} \max_{(\bsig,\btau,\bmu) \in \calZ_a^+}
\{\half \by^T \BB \by -\bh^T \by +(\bz^*)^T \bz\} \label{minmaxza}
\end{eqnarray*}
Note that
\begin{eqnarray*}
\max_{\bz^* \in \calZ_a^+} \{ \WW^{\sharp}(\bz^*)\}=0
\end{eqnarray*}
and
\begin{eqnarray*}
\max_{\bz\in \ZZ} \{ \WW(\bz)\}=0.
\end{eqnarray*}
Thus, it follows from (\ref{minmaxza}) that
\begin{eqnarray*}
P^d(\barbsig,\barbtau,\barbmu)&=&
\min_{\by \in \calY} \max_{(\bsig,\btau, \bmu)\in \calZ_a^+}
\{\half \by^T \BB \by -\bh^T \by+(\bz^*)^T \bz -\WW^{\sharp} (\bz^*)\}\\
&=& \min_{\by \in \calY} \{\half \by^T \BB \by -\bh^T \by \}
+ \max_{(\bsig,\btau,\bu) \in \calZ_a^+ }
\{ (\bz^*)^T \bz- \WW^{\sharp} (\bz^*)\}\nonumber\\
&=& \min_{\by \in \calY} \{\half \by^T \BB \by -\bh^T \by\}
+ \max_{(\bsig,\btau,\bmu) \in \calZ_a^+}
\{(\bz^*)^T \bz -\bz^T \bz^* +\WW(\bz)\}\nonumber\\
&=& \min_{\by \in \calY}\{ \half \by^T \BB \by -\bh^T \by \}\nonumber\\
&=& \min_{\by \in \calY} \PP(\by).
\end{eqnarray*}
This completes the proof.
\hfill $\Box$
\section{Numerical Experience}
All data and computational results presented in this section are produced
by Matlab. In order to save space and fit the matrix
in the paper, we round our these results  up to two decimals.

{\bf Example 1. 5-dimensional problem}.

Consider Problem $(\calP_a)$ with
$\bx=[x_1,\cdots,x_5]^T$ , while $x_i \in \{ 2,3,5 \} $, $i=1, \cdots, 5$,
\begin{eqnarray*}
\bQ=\left[
\begin{array}{ccccc}
3.43&0.60&0.39&0.10&0.60\\
0.60&2.76&0.32&0.65&0.49\\
0.39&0.32&2.07&0.59&0.39\\
0.10&0.65&0.59&2.62&0.30\\
0.60&0.49&0.39&0.30&3.34
\end{array}
\right],
\end{eqnarray*}
\begin{eqnarray*}
\bc=[38.97,-24.17,40.39,-9.65,13.20]^T,
\end{eqnarray*}
\begin{eqnarray*}
\bA=\left[
\begin{array}{ccccc}
0.94&0.23&0.04&0.65&0.74\\
0.96&0.35&0.17&0.45&0.19\\
0.58&0.82&0.65&0.55&0.69\\
0.06&0.02&0.73&0.30&0.18
\end{array}
\right],
\end{eqnarray*}
\begin{eqnarray*}
\bb=[11.49,9.32,14.43,5.66]^T.
\end{eqnarray*}
Under the transformation (\ref{01transfer}), this problem
is transformed into the  0-1 programming Problem $(\calP)$, where
\begin{eqnarray*}
\by=[y_{1,1}, y_{1,2}, y_{1, 3},\cdots, y_{5,1}, y_{5,1},y_{5,3}]^T
\in \real^{15},
\end{eqnarray*}
\scriptsize
\begin{eqnarray*}
\bB=
\left[
\begin{array}{ccccccccccccccc}
13.71&20.56&34.27&2.40 &3.61 &6.01 &1.58  &2.37 &3.95 &0.39 &0.58 &0.97 &2.38 &3.57 &5.95\\
20.56&30.84&51.41&3.61 &5.41 &9.01 &2.37  &3.55 &5.92 &0.58 &0.88 &1.46 &3.57 &5.36 &8.93\\
34.27&51.41&85.68&6.01 &9.01 &15.02&3.95  &5.92 &9.87 &0.97 &1.46 &2.43 &5,95 &8.93 &14.88\\
2.40 &3.61 &6.01 &11.05&16.57&27.61&1.27  &1.91 &3.18 &2.61 &3.91 &6.52 &1.95 &2.93 &4.88\\
3.61 &5.41 &9.01 &16.57&24.85&41.42&1.91  &2.86 &4.77 &3.91 &5.87 &9.78 &2.93 &4.39 &7.32\\
6.01 &9.01 &15.02&27.61&41.42&69.03&3.18  &4.77 &7.96 &6.52 &9.78 &16.31&4.88 &7.32 &12.20\\
1.58 &2.37 &3.95 &1.27 &1.91 &3.18 &8.27  &12.40&20.67&2.37 &3.55 &5.92 &1.57 &2.36 &3.93\\
2.37 &3.55 &5.92 &1.91 &2.86 &4.77 &12.40 &18.60&31.00&3.55 &5.33 &8.89 &2.36 &3.53 &5.90\\
3.95 &5.92 &9.87 &3.18 &4.77 &7.96 &20.67 &31.00&51.67&5.92 &8.86 &14.81&3.93 &5.90 &9.83\\
0.39 &5.58 &0.97 &2.61 &3.91 &6.52 &2.37  &3.55 &5.92 &10.50&15.74&26.24&1.20 &1.80 &3.00\\
0.58 &0.88 &1.46 &3.91 &5.87 &9.78 &3.55  &5.33 &8.89 &15.74&23.62&39.36&1.80 &2.70 &4.50\\
0.97 &1.46 &2.43 &6.52 &9.78 &16.31&5.92  &8.89 &14.81&26.24&39.36&65.60&3.00 &4.50 &7.51\\
2.38 &3.57 &5.95 &1.95 &2.93 &4.88 &1.57  &2.36 &3.93 &1.20 &1.80 &3.00 &13.35&20.02&33.37\\
3.57 &5.36 &8.93 &2.93 &4.39 &7.32 &2.36  &3.54 &5.90 &1.80 &2.70 &4.50 &20.02&30.04&50.06\\
5.95 &8.93 &14.88&4.88 &7.32 &12.20&3.93  &5.90 &9.83 &3.00 &4.50 &7.51 &33.37&50.06&83.43
\end{array}
\right],
\end{eqnarray*}
\normalsize
\begin{eqnarray*}
\bh &=& [77.95,116.92,194.87,-48.34,-72.51,-120.85,80.78,121.17\\
&& 201.96,-19.29,-28.94, -48.23,26.39,39.59,65.99]^T,
\end{eqnarray*}
\scriptsize
\begin{eqnarray*}
\bD=\left[
\begin{array}{ccccccccccccccc}
\;1.88\;\;2.83\;\;4.71\;\;0.47\;\;0.70\;\;1.17\;\;0.09\;\;0.12\;\;0.22\;\;1.30\;\;1.94\;\;3.24\;\;1.49\;\;2.23\;\;3.72\\
\;1.91\;\;2.87\;\;4.78\;\;0.71\;\;1.06\;\;1.77\;\;0.34\;\;0.51\;\;0.85\;\;0.90\;\;1.35\;\;2.25\;\;0.38\;\;0.57\;\;0.94\\
\;1.15\;\;1.72\;\;2.88\;\;1.64\;\;2.46\;\;4.11\;\;1.30\;\;1.95\;\;3.25\;\;1.09\;\;1.64\;\;2.74\;\;1.37\;\;2.06\;\;3.43\\
\;0.12\;\;0.18\;\;0.30\;\;0.03\;\;0.05\;\;0.08\;\;1.46\;\;2.20\;\;3.66\;\;0.59\;\;0.89\;\;1.48\;\;0.37\;\;0.55\;\;0.92
\end{array}
\right],
\end{eqnarray*}
\normalsize
\begin{eqnarray*}
\bH =
\left[
\begin{array}{cccccccccc}
1  & \cdots & 1 & 0 &\cdots& 0 &\cdots& 0 &\cdots &0\\
0  & \cdots & 0 & 1&\cdots& 1 &\cdots& 0 &\cdots &0\\
\vdots&\ddots&\vdots&\vdots&\ddots&\vdots&\ddots&\vdots&\ddots&\vdots\\
0  & \cdots & 0 & 0 &\cdots& 0 &\cdots& 1 &\cdots &1
\end{array}
\right]\in \real^{5 \times 15}.
\end{eqnarray*}
The Canonical dual problem can be stated as follows:
\begin{eqnarray*}
&(\calP^{d})\;\;&{\rm Maximize}\;\;P^d(\bsig, \btau,\bmu)
=-\half F(\bsig, \btau, \bmu)^T G^\dagger(\bmu) F(\bsig, \btau, \bmu)
-\bsig^T \bb-\btau^T \be_{5} \\
&& {\rm subject \; to}\;\; \bsig \ge 0, \bmu >0,\\
&&\;\;\;\;\;\;\;\;\;\;\;\;\;\;\;\;\;\bsig \in \real^4, \btau \in \real^5,
\bmu \in \real^{15},\nonumber
\end{eqnarray*}
where $F(\bsig, \btau, \bmu)$ and $G(\bmu)$ are as defined by
(\ref{defg}) and (\ref{deff}), respectively.\\
By solving this dual problem with the sequential quadratic
programming method in the optimization Toolbox within the
Matlab environment, we obtain
\begin{eqnarray*}
\barbsig=[0,0,0,0]^T,
\end{eqnarray*}
\begin{eqnarray*}
\barbtau=[73.90,-106.70,111.95,-59.27,-0.01]^T,
\end{eqnarray*}
and
\begin{eqnarray*}
\barbmu &=& [39.34,22.07,12.49,33.56,3.01,76.14,61.00,35.52\\
&& 18.78,1.47,41.96, 0.001,0.001,0.006]^T.
\end{eqnarray*}
It is clear that
 $(\barbsig,\barbtau,\barbmu) \in \calZ_a^+$. Thus, from Theorem \ref{pdthm},
\begin{eqnarray*}
\barby&=& (\BB +2 \Diag (\barbmu))^\dagger(\bh - \DD^T \barbsig -\HH^T \barbtau + \barbmu)\\
&=&[0,0,1,1,0,0,0,0,1,1,0,0,1,0,0]^T
\end{eqnarray*}
is the global minimizer of Problem $(\calP)$ with
$\PP^d(\barbsig, \barbtau,\barbmu)=-227.87=\PP(\barby)$.
The solution to the original primal problem
can be calculated by using the transformation
\begin{eqnarray*}
\barx_i=\sum_{j=1}^{K_i} u_{i,j} \bary_{i,j},\; i=1,2,3,4,5,
\end{eqnarray*}
to give
\begin{eqnarray*}
\barbx=[5,2,5,2,2]^T
\end{eqnarray*}
with $P(\barbx)=-227.87$.

{\bf Example 2. 10-dimensional problem}.

Consider Problem $(\calP_a)$, with $\bx=[x_1, \cdots, x_{10}]^T $, while
$x_i \in \{1,2,4,7,9\},\;i=1,\cdots,10$,
\begin{eqnarray*}
\bQ=\left[
\begin{array}{cccccccccc}
6.17&0.62&0.46&0.37&0.56&0.66&0.67&0.85&0.57&0.44\\
0.62&5.63&0.29&0.56&0.79&0.29&0.43&0.69&0.49&0.39\\
0.46&0.29&5.81&0.55&0.22&0.55&0.36&0.27&0.51&0.91\\
0.37&0.56&0.55&6.10&0.28&0.42&0.44&0.34&0.75&0.44\\
0.56&0.79&0.22&0.28&4.75&0.40&0.55&0.42&0.49&0.44\\
0.66&0.29&0.55&0.42&0.40&5.71&0.32&0.57&0.65&0.70\\
0.67&0.43&0.36&0.44&0.55&0.32&5.27&0.56&0.37&0.85\\
0.85&0.69&0.27&0.34&0.42&0.57&0.56&5.91&0.15&0.62\\
0.57&0.49&0.51&0.75&0.49&0.65&0.37&0.15&4.51&0.46\\
0.44&0.39&0.91&0.44&0.44&0.70&0.85&0.62&0.46&5.73
\end{array}
\right],
\end{eqnarray*}
\begin{eqnarray*}
\bff=[0.89,0.03, 0.49, 0.17, 0.98, 0.71, 0.50, 0.47, 0.06, 0.68]^T,
\end{eqnarray*}
\begin{eqnarray*}
\bA=\left[
\begin{array}{cccccccccc}
0.04&0.82&0.97&0.83&0.83&0.42&0.02&0.20&0.05&0.94\\
0.07&0.72&0.65&0.08&0.80&0.66&0.98&0.49&0.74&0.42\\
0.52&0.15&0.80&0.13&0.06&0.63&0.17&0.34&0.27&0.98\\
0.10&0.66&0.45&0.17&0.40&0.29&0.11&0.95&0.42&0.30\\
0.82&0.52&0.43&0.39&0.53&0.43&0.37&0.92&0.55&0.70
\end{array}
\right],
\end{eqnarray*}
\begin{eqnarray*}
\bb=[33.76, 37.07, 26.75, 25.46, 37.36]^T.
\end{eqnarray*}
By solving the canonical dual problem of Problem $(\calP_a)$,
 we obtain
\begin{eqnarray*}
\barbsig=[0,0,0,0,0]^T,
\end{eqnarray*}
\begin{eqnarray*}
\barbtau &=& [-19.99, -20.12, -18.13, -18.37, -14.32,\\
&& -17.13, -18.46, -19.73,
-17.65, -16.55]^T,
\end{eqnarray*}
and
\begin{eqnarray*}
\barbmu &=&
[9.51,0.97,21.93,53.36,74.34,9.95,0.21,20.53,51.01,71.35\\
&& 8.68,0.77,19.68,48.03,66.94,8.30,1.77,21.91,52.13,72.27\\
&& 6.40,1.54,17.39,41.19,57.04,7.57,1.98,21.10,49.77,68.90\\
&& 9.15,0.16,18.79,46.72,65.34,9.82,0.09,19.90,49.63,69.45\\
&& 8.76,0.13,17.92,44.60,62.39,6.26,4.03,24.60,55.48,76.04]^T,
\end{eqnarray*}
It is clear that $(\barbsig,\barbtau,\barbmu) \in \calZ_a^+$. Therefore,
\begin{eqnarray*}
&&\barby=[1,0,0,0,0,1,0,0,0,0,1,0,0,0,0,1,0,0,0,0,1,0,0,0,0,\\
&&1,0,0,0,0,1,0,0,0,0,1,0,0,0,0,1,0,0,0,0,1,0,0,0,0]^T
\end{eqnarray*}
is the global minimizer of the problem $(\calP)$ with
 $\PP^d(\barbsig, \barbtau,\barbmu)=45.54=\PP(\barby)$.
The solution to the original primal problem is
\begin{eqnarray*}
\barbx=[1,1,1,1,1,1,1,1,1,1]^T
\end{eqnarray*}
with $P(\barbx)=45.54$.

{\bf Example 3. Large scale problems}.

Consider Problem $(\calP_a)$ with $n=20$, $50$, $100$, $200$ and $300$. Let
these five problems be referred to as Problem (1), $\cdots$, Problem (5),
respectively. Their coefficients are generated randomly with uniform distribution.
For each problem, $q_{ij} \in (0,1)$, $ a_{ij} \in (0,1)$, for
$i=1, \cdots, n$; $j=1, \cdots, n$, and  $c_i \in (0,1)$, $x_i \in \{1,2,3,4,5 \}$,
for $i=1, \cdots n$. Without loss of generality, we
ensure that the constructed  $\QQ $ is a symmetric matrix.
Otherwise, we let $\QQ= \frac{\QQ +\QQ^T }{2}$. Furthermore, let
$\QQ$ be such that it is diagonally dominated.
For each $x_i$, its lower bound is $l_i=1$, and its upper
bound is $u_i=5$. Let $l=[l_1, \cdots, l_n]^T$ and  $u=[u_1,\cdots,u_n]^T $.
The right-hand sides of the linear constraints are chosen such that the
feasibility of the test problem is satisfied. More specifically, we
 set $\bb=\sum_{j} a_{ij} l_j +
0.5\cdot (\sum_j a_{ij}u_j-\sum_j a_{ij} l_j)$.

We then construct the canonical problem of each of the five problems.
It is solved by using the sequential quadratic programming method with
active set strategy from the Optimization Toolbox within the Matlab environment.
The specifications of the personal notebook computer used are:
Intel(R), Core(TM)(1.20 GHZ), Window Vista(TM).
Table 1 presents the numerical results,
where $m$ is number of linear constraints in Problem $I(\calP_a)$.
\begin{table}[h!]
\vspace*{0.1in}
\caption{Numerical results for large scale integer programming problems }
\centering
\begin{tabular}{cccc}
\hline
n& m& CPU Time (Seconds)\\
\hline
20&5&4.8\\
50&5&19.1\\
100&5&75.8\\
200&5&277.8\\
300&5&649.7\\
\hline
\end{tabular}
\end{table}

From Table 1, we see that the algorithm based on the
canonical dual method can solve
large scale problems with reasonable computational time.
Furthermore, for each of the five problems, the solution obtained
is a global optimal solution.
For the
case of $n=300$, the
equivalent problem
in the form of Problem $(\calP_b)$ has  1500
variables. For  such a problem,  there are  $2^{1500}$
possible  combinations.

\section{Conclusion}
We have presented a canonical duality theory for solving a general
discrete value selection problem with quadratic cost function
and linear constraints. Our results show that this
NP-hard problem can be converted to a continuous
concave dual maximization problem  without duality gap.
If the canonical dual space $\calZ^+_a$ is non empty,
the problem can be solved easily via well-developed convex optimization methods.
Several examples, including some large scale ones, were solved
effectively by using the method proposed.

\vspace{1.0 cm}

\noindent{\bf Acknowledgement}:
This paper was partially supported by a grant (AFOSR FA9550-10-1-0487)
from the US Air Force Office of Scientific Research. Dr. Ning Ruan was
supported by a funding from the Australian Government
under the Collaborative Research Networks (CRN) program.

\end{document}